\documentstyle[12pt]{article}
\newcommand{\res}{\mathop{\rm res}\nolimits}

\newcommand{\bg}{\phantom._0\Phi_1(-;0;q^2,i(1-q^2)q^2zs)}
\newcommand{\bgm}{\phantom._0\Phi_1(-;0;q^2,-i(1-q^2)q^2zs)}
\newcommand{\Ex}{E_{q^2}(-i(1-q^2)zs)}
\newcommand{\mEx}{E_{q^2}(i(1-q^2)zs)}

\newcommand{\La}{\Lambda}
\newcommand{\G}{\Gamma}
\newcommand{\p}{\partial}

\newtheorem{predl}{Proposition}[section]
\newtheorem{defi}{Definition}[section]

\newtheorem{cor}{Corollary}[section]

\newtheorem{theo}{Theorem}[section]
\newcommand{\beq}[1]{\begin{equation}\label{#1}}
\newcommand{\eq}{\end{equation}}

\begin{document}

\vspace{10mm}
\begin{flushright}
 ITEP-TH-/00\\
\end{flushright}
\vspace{10mm}
\begin{center}
{\Large \bf $q$-convolution and its $q$-Fourier transform}\\
\vspace{5mm}
V.-B.K.Rogov \footnote{The author was supported in part by the
Russian Foundation for Fundamental Research (grant no. 0001-00143) and
the NIOKR MPS RF} \\
MIIT, 101475, Moscow\\ e-mail vrogov@cemi.rssi.ru\\

\end{center}
\begin{abstract}
The functions on the lattice generated by the integer degrees of $q^2$
are considered, $0<q<1$. The $q^2$-translation operator is defined.
The multiplicators and the $q^2$-convolutors are defined in the functional
spaces which are dual with respect to the $q^2$-Fourier transform. The
$q^2$-analog of convolution of two $q^2$-distributions is constructed.
The $q^2$-analog of an arbitrary order derivative is introduced
\end{abstract}

\section{Introduction}
\setcounter{equation}{0}

The Fourier transform plays an important part in the harmonic analysis
on the simple Lie groups and on the homogeneous spaces. The concept of the
convolution is closely connected with the Fourier transform because the
last moves the convolution of two functions to the product of
their images. In case of the quantum groups and quantum homogeneous
spaces the $q$-analogs of the Fourier transform play the same part and the
problem to construct the $q$-analog of the convolution arises.
The different $q$-analogs of the Fourier transform have been investigated
in \cite{CZ,AAS,KM,KS}. The $q$-convolution was introduced for the first
time in \cite{KM}. It is extensive investigated in \cite{CK,C}. The
definition of the $q$-convolution is connected with the definition of the
$q$-Fourier transform because the $q$-Fourier transform moves the
$q$-convolution of two functions to the product of their images. The
$q$-convolution considered in \cite{CK,C} is connected with the
$q$-Fourier transform considered in \cite{KS}. In these works the braided
line is introduced.

In \cite{OR} the $q^2$-Fourier transform and the inversion formula have
been constructed and they are quite similar to the classical ones
\cite{GS1,GS2}. This construction coincides with the classical
Fourier transform if $q\to1$. The construction of the $q^2$-convolution
corresponding the $q^2$-Fourier transform in  the space of
$q^2$-distributions \cite{OR} is proposed in the present work.  Thus
constructing of the theory similar to the classical one \cite{GS1,GS2} is
prolonged. Moreover the braided line construction presents implicitly in
this work, because we consider non commuting variables.

The $q^2$-convolution operator allows to determine the $q^2$-derivative of
an arbitrary order.

In this paper we will use the same notation that in \cite{OR}.

\section{Some preliminary relations}
\setcounter{equation}{0}

We assume that $z\in{\bf C}$ and $|q|<1$, unless otherwise is specified.

We recall some notations \cite{GR}. For an arbitrary $a$
$$
(a,q)_n=\left\{
\begin{array}{lcl}
1&{\rm for}&n=0\\
(1-a)(1-aq)\ldots(1-aq^{n-1})&{\rm for}&n\ge1,\\
\end{array}
\right.
$$
$$
(a,q)_\infty=\lim_{n\to\infty}(a,q)_n,
$$
$$
\left[\begin{array}{c}
l \\i
\end{array}
\right ]_{q^2}
=\frac
{(q^2,q^2)_l}
{(q^2,q^2)_i(q^2,q^2)_{l-i}}.
$$
Consider the $q^2$-exponentials
\begin{equation}\label{2.1}
e_{q^2}(z)=\sum_{n=0}^\infty\frac{z^n}{(q^2,q^2)_n}=
\frac1{(z,q^2)_\infty},\qquad|z|<1,
\end{equation}
\begin{equation}\label{2.2}
E_{q^2}(z)=\sum_{n=0}^\infty\frac{q^{n(n-1)}z^n}{(q^2,q^2)_n}=
(-z,q^2)_\infty
\end{equation}
and the basic hypergeometric series
$$
\phantom._r\Phi_s(a_1,\ldots,a_r;b_1,\ldots,b_s;q^2,z)=
$$
\begin{equation}\label{2.3}
=\sum_{n=0}^\infty\frac{(a_1,q^2)_n\cdots(a_r,q^2)_n}
{(q^2,q^2)_n(b_1,q^2)_n\cdots(b_s,q^2)_n}
\lbrack(-1)^nq^{n\choose2}\rbrack^{1+s-r} z^n.
\end{equation}

We consider the series
$$
{\bf Q}(z,q)=(1-q^2)\sum_{m=-\infty}^\infty\frac1
{zq^{2m}+z^{-1}q^{-2m}}
$$
expressed by theta-function (see \cite{OR}). Assume
\begin{equation}\label{2.4}
\Theta_0={\bf Q}(1-q^2,q).
\end{equation}

Let ${\cal A}=C(z,z^{-1})$ be the algebra of formal Laurent series. The
$q^2$-derivative of $f(z)\in {\cal A}$ is defined as
$$
\p_zf(z)=\frac{z^{-1}}{1-q^2}(f(z)-f(q^2z)).
$$
For an arbitrary $n\ge0$
\begin{equation}\label{2.5}
\p_z^kz^n=\left\{
\begin{array}{lcl}
\frac{(q^2,q^2)_n}{(q^2,q^2)_{n-k}(1-q^2)^k}z^{n-k}&{\rm for}&0\le k\le n\\
0&{\rm for}&k>n,\\
\end{array}
\right.
\end{equation}
and for any $n\ge0$ and $k\ge0$
\begin{equation}\label{2.6}
\p_z^kz^{-n-1}=(-1)^kq^{-k(2n+k+1)}\frac{(q^2,q^2)_{n+k}}
{(q^2,q^2)_n(1-q^2)^k}z^{-n-k-1}.
\end{equation}

The $q^2$-integral (Jackson integral \cite{GR}) is defined as the map
$I_{q^2}$ from ${\cal A}$ to the space of formal number series
$$
I_{q^2}f=\int d_{q^2}zf(z)=(1-q^2)\sum_{m=-\infty}^\infty q^{2m}
[f(q^{2m})+f(-q^{2m})]
$$
\begin{defi}\label{d2.1}
$f(z)$ is absolutely $q^2$-integrable function, if the series
$$
\sum_{m=-\infty}^\infty q^{2m}[|f(q^{2m})|+|f(-q^{2m})|]
$$
converges.
\end{defi}

Let ${\cal B}$ be the algebra analogous to ${\cal A}$, but generated by
$s,s^{-1}$ which commute with $z$ as $zs=q^2sz$.

We denote by ${\cal AB}$ the whole algebra with generators
$z,z^{-1},s,s^{-1}$ and relations
\begin{equation}\label{2.7}
zs=q^2sz,\quad \p_zs=q^{-2}s\p_z,\quad \p_sz=q^2z\p_s,
\quad \p_z\p_s=q^2\p_s\p_z.
\end{equation}

We will consider ${\cal AB}$ as a left module under the action of
${\cal A}$ by multiplication, and a right module under the action of
${\cal B}$.

To define the $q^2$-integral on ${\cal AB}$ we order the generators of
integrand in such a way that $z$ stays on the left side while $s$ stays on
the right side. For example, if $f(z)=\sum_ra_rz^r$ then
$$
f(zs)=\sum_ra_r(zs)^r=\sum_ra_rq^{-r(r-1)}z^rs^r.
$$
For convenience we introduce the following notation
$$
\ddag g(zs)\ddag=\sum_ra_rz^rs^r,~{\rm if}~g(z)=\sum_ra_rz^r.
$$
For example, we can derive from (\ref{2.1}),  (\ref{2.2}) and (\ref{2.7})
that
$$
E_{q^2}((1-q^2)zs)=\ddag e_{q^2}((1-q^2)zs)\ddag.
$$

We define the operator
$$
\Lambda_z: f(z) \to f(q^2z).
$$
Obviously
$$
\La_zz=q^2z\La_z,\qquad \p_z\La_z=q^2\La_z\p_z.
$$

\section{$q^2$-distributions and the $q^2$-Fourier transform}
\setcounter{equation}{0}

In \cite{OR} we have defined the spaces of the test functions and
the $q^2$-distributions and the $q^2$-Fourier transform has been
constructed.  We reproduce some statements here.

Let $S_{q^2}=\{\phi(x)\}$ be the space of infinitely $q^2$-differentiable
fast decreasing functions
\begin{equation}\label{3.1}
|x^k\p_x^l\phi(x)|\le C_{k,l}(q),~k\ge0,~ l\ge0.
\end{equation}

Let $S$ be the space of infinitely differentiable (in the classic sense)
fast decreasing function
$$
|x^k\phi^{(l)}(x)|\le C_{k,l},~k\ge0~l\ge0.
$$
It has been show in \cite{OR} that $S\subset S_{q^2}$. In addition
\begin{equation}\label{3.2}
\p_x^l\phi(x)|_{x=0}=\frac{(q^2,q^2)_l}{(1-q^2)^ll!}\phi^{(l)}(0).
\end{equation}

\begin{predl}\label{p3.1}
If $\phi(z)\in S_{q^2}$, then
$$
\int d_{q^2}z\p_z\phi(z)=0.
$$
\end{predl}
\begin{cor}\label{c3.1}
($q^2$-integration by parts). For any $k\ge 0$
\begin{equation}\label{3.3}
\int d_{q^2}z\phi(z)\p_z^k\psi(z)=
(-1)^kq^{-k(k-1)} \int d_{q^2}z\p_z^k\phi(z)\psi(q^{2k}z).
\end{equation}
\end{cor}

\begin{defi}\label{d3.1}
The $q^2$-distribution $f$ over $S_{q^2}$ is a linear continuous
functional
$$
f \to <f,\phi>, ~~\phi(z)\in S_{q^2}.
$$
\end{defi}

We denote by $S_{q^2}'$ the space of the $q^2$-distributions over
$S_{q^2}$.

The $q^2$-distributions defined by the $q^2$-integral
$$
<f,\phi>=\int_{-\infty}^\infty d_{q^2}z\overline{f(z)}\phi(z)=
(1-q^2)\sum_{m=-\infty}^\infty q^{2m}
[\overline{f(q^{2m})}\phi(q^{2m})+\overline{f(-q^{2m})}\phi(-q^{2m})]
$$
we refer as a regular one.

Proposition \ref{p3.1} and Corollary \ref{c3.1} allow to introduce the
$q^2$-differentiation in $S_{q^2}'$
\begin{equation}\label{3.4}
<\p_zf,\phi>=-<\La_z f,\p_z\phi>.
\end{equation}

It follows from (\ref{3.4}) that the conjugate operator for $\p_z^k$
for any $k\ge0$ has the form
\begin{equation}\label{3.5}
(\p_z^k)^*=(-1)^kq^{k(k-1)}\p_z^k\La_z^{-k}.
\end{equation}
The change of variables $q^{-2k}z\to z$ in the $q^2$-integral leads to
\begin{equation}\label{3.6}
(\La_z^{-k})^*=q^{2k}\La_z^k.
\end{equation}

\begin{defi}\label{d3.2}
$f$ is the $q^2$-distribution with multiplicity $p$ of the
$q^2$-singularity if it is represented in form
$$
f=\sum_{k=0}^p\p_z^kf_k(z),
$$
where $f_k(z)$ are the ordinary functions growing no faster then some
power of $|z|$ as $|z|\to\infty$.
\end{defi}

For example $\delta_{q^2}(z)$ is the $q^2$-distribution of
multiplicity one of the $q^2$-singularity because for an arbitrary
$\phi(z)\in S_{q^2}$

$$
<\delta_{q^2},\phi>=<\frac12\p_z(\theta_{q^2}^+-\theta_{q^2}^-),\phi>=
-<\frac12\Lambda_z(\theta_{q^2}^+-\theta_{q^2}^-),\p_z\phi>=
=-\frac12\int_0^\infty d_{q^2}z\p_z\phi(z)+
$$
$$
+\frac12\int_{-\infty}^0d_{q^2}z\p_z\phi(z)=
-\frac12\sum_{m=-\infty}^\infty[\phi(q^{2m})-\phi(q^{2m+2})+
\phi(-q^{2m})-\phi(-q^{2m+2})]=
$$
$$
=\lim_{m\to\infty}\frac{\phi(q^{2m})+\phi(-q^{2m})}2=\phi(0).
$$

Let the space $S^{q^2}=\{\psi(s)\}$ be the copy of the
$S_{q^2}=\{\phi(z)\}$ (\ref{3.1}), but $s$ and $z$ behave
as the generators of the algebra ${\cal AB}$ (\ref{2.7}). Introduce the
same topology in $S^{q^2}$ as one in $S_{q^2}$
$$
|s^k\p_s^l\phi(s)|\le C_{k,l}(q),~k\ge0,~ l\ge0 ,
$$
Thereby these spaces are isomorphic.

The $q^2$-Fourier transform ${\cal F}_{q^2}$, i.e. the map $S_{q^2}$ into
$S^{q^2}$ has been constructed in \cite{OR}
$$
S_{q^2} \stackrel{{\cal F}_{q^2}}{\longrightarrow} S^{q^2}
$$
where
$$
{\cal F}_{q^2}\phi(z)=\int d_{q^2}z\phi(z)\bg,
$$
and $\phantom._0\Phi_1$ is determined by (\ref{2.3}).
The inverse transform
\begin{equation}\label{3.7}
{\cal F}_{q^2}^{-1}\psi(s)=\frac1{2\Theta_0}\int E_{q^2}(-i(1-q^2)zs)
\psi(s)d_{q^2}s, ~~~\psi(s)\in S^{q^2},
\end{equation}
has been constructed and their continuity was proofed. The constant
$\Theta_0$ is determined by (\ref{2.4}).

The following relations are valid
$$
{\cal F}_{q^2}\La_z=q^{-2}\La_s^{-1}{\cal F}_{q^2},
\qquad {\cal F}_{q^2}\p_z=-is{\cal F}_{q^2},
\qquad {\cal F}_{q^2}z=-iq^{-2}\La_s^{-1}\p_s{\cal F}_{q^2},
$$
\begin{equation}\label{3.8}
{\cal F}_{q^2}^{-1}\La_s=q^{-2}\La_z^{-1}{\cal F}_{q^2}^{-1},
\qquad {\cal F}_{q^2}^{-1}\p_s=i\La_z^{-1}z{\cal F}_{q^2}^{-1},
\qquad {\cal F}_{q^2}^{-1}s=i\p_z{\cal F}_{q^2}^{-1}.
\end{equation}
\begin{defi}\label{d3.3}
The $q^2$-Fourier transform of a $q^2$-distribution $f\in S_{q^2}'$ is
the $q^2$-distribution $g\in(S^{q^2})'$ defined by the equality
\begin{equation}\label{3.9}
<g,\psi>=<f,\phi>, ~~~\psi(s)={\cal F}_{q^2}\phi(z),
\end{equation}
where $\phi(z)$ is an arbitrary function from $S_{q^2}$.
\end{defi}
Presuppose that the $q^2$-distribution $f$ corresponds to $f(z)$ and
$zf(z)$ is absolutely $q^2$-integrable function. Let
$\phi(z)= {\cal F}_{q^2}^{-1}\psi(s)$. Then
$$
<f,\phi>=\frac1{2\Theta_0}
\int d_{q^2}z\overline{f(z)}\int\Ex\psi(s)d_{q^2}s=
$$
$$
=\frac1{2\Theta_0}\int\overline{\int d_{q^2}zf(z)\mEx}
\psi(s)d_{q^2}s=<g,\psi>.
$$
It means that the $q^2$-distribution $g$ corresponds to the function
\begin{equation}\label{3.10}
g(s)=\frac1{2\Theta_0}\int d_{q^2}zf(z)\mEx.
\end{equation}

In the same way, if $g$ is determined by the absolutely
$q^2$-integrable function $g(s)$ and $\psi(s)={\cal F}_{q^2}\phi(z)$, then
$$
<g,\psi>=\int\int d_{q^2}z\phi(z)\bg\overline{g(s)}d_{q^2}s=
$$
$$
=\int d_{q^2}z\phi(z)\overline{\int\bgm g(s)d_{q^2}s}=<f,\phi>,
$$
i.e. $f$ corresponds to
\begin{equation}\label{3.11}
f(z)=\int\bgm g(s)d_{q^2}s.
\end{equation}

The $q^2$-Fourier transform of a $q^2$-distribution from $S_{q^2}'$
we denote by ${\cal F}_{q^2}'$. It follows from (\ref{3.10}), (\ref{3.11})
and (\ref{3.8}) that in the space of $q^2$-distributions the following
commutative relations are valid
$$
{\cal F}_{q^2}'\La_z=q^{-2}\La_s^{-1}{\cal F}_{q^2}', \qquad
{\cal F}_{q^2}'\p_z=-i\La_s^{-1}s{\cal F}_{q^2}', \qquad
{\cal F}_{q^2}'z=-i\p_s{\cal F}_{q^2}',
$$
\begin{equation}\label{3.12}
({\cal F}_{q^2}')^{-1}\La_s=q^{-2}\La_z^{-1}({\cal F}_{q^2}')^{-1},
({\cal F}_{q^2}')^{-1}\p_s=iz({\cal F}_{q^2}')^{-1},
({\cal  F}_{q^2}')^{-1}s= iq^{-2}\La_z^{-1}\p_z({\cal F}_{q^2}')^{-1}.
\end{equation}

\section{$q^2$-shift in the space of the test functions}
\setcounter{equation}{0}

Let $\xi$ be an element of the same nature as $s$ so their sum is
determined
$$
s+\xi=\xi+s,
$$
and the commutative relations
\begin{equation}\label{4.1}
\xi s=q^2s\xi, \quad \xi\p_s=q^{-2}\p_s\xi,\quad \La_s\xi=\xi\La_s
\end{equation}
are fulfilled. In this case we will call the element $s$ subordinate to
$\xi$.
\begin{defi}\label{d4.1}
We will call the operator
\begin{equation}\label{4.2}
T_\xi=e_{q^2}((1-q^2)\xi\La_s^{-1}\p_s)
\end{equation}
by $q^2$-shift in the space $S^{q^2}$.
\end{defi}
\begin{predl}\label{p4.1}
For an arbitrary function $g(s)\in\cal B$
$$
T_\xi g(s)=g(s+\xi).
$$
\end{predl}
{\bf Proof.} It follows from (\ref{2.1}), (\ref{4.1}) and
(\ref{2.6}) that for an arbitrary $n\ge0$
\begin{equation}\label{4.3}
e_{q^2}((1-q^2)\xi\La_s^{-1}\p_s)s^n=
\ddag e_{q^2}((1-q^2)\xi\La_s^{-1}\p_s)\ddag s^n=\sum_{k=0}^n
\left[\begin{array}{c}
n \\k
\end{array}
\right ]_{q^2}
s^{n-k}\xi^k=(s+\xi)^n.
\end{equation}
By induction on $n$ for any $n\ge0$ one finds
$$
(s+\xi)^{-n-1}=\sum_{k=0}^\infty(-1)^k
\left[\begin{array}{c}
n+k \\k
\end{array}
\right ]_{q^2}
s^{-n-k-1}\xi^k.
$$
Therefore, it follows from (\ref{2.1}), (\ref{4.1}) and (\ref{2.7})
for any $n\ge0$
\begin{equation}\label{4.4}
e_{q^2}((1-q^2)\xi\La_s^{-1}\p_s)s^{-n-1}
=\ddag e_{q^2}((1-q^2)\xi\La_s^{-1}\p_s)\ddag s^{-n-1}=(s+\xi)^{-n-1}.
\end{equation}
Now it follows from (\ref{4.3}) and (\ref{4.4}) that for
$g(s)=\sum_ra_rs^r$
$$
e_{q^2}((1-q^2)\xi\La_s^{-1}\p_s)g(s)=g(s+\xi).
$$
\rule{5pt}{5pt}
\begin{predl}\label{p4.2}
If $\xi_2$ is subordinated $\xi_1$, i.e.
$$
\xi_1\xi_2=q^2\xi_2\xi_1,
$$
$$
T_{\xi_2}T_{\xi_1}=T_{\xi_1+\xi_2}.
$$
\end{predl}
{\bf Proof}. It follows from (\ref{2.1}) and (\ref{5.1})
$$
e_{q^2}((1-q^2)\xi_2\La_s^{-1}\p_s)e_{q^2}((1-q^2)\xi_1\La_s^{-1}\p_s)=
\sum_{k=0}^\infty\frac{(1-q^2)^k}{(q^2,q^2)_k}\xi_2^k\La_s^{-k}\p_s^k
\sum_{l=0}^\infty\frac{(1-q^2)^l}{(q^2,q^2)_l}\xi_1^l\La_s^{-l}\p_s^l=
$$
$$
=\sum_{m=0}^\infty\frac{(1-q^2)^m}{(q^2,q^2)_m}
\sum_{l=0}^m
\left[\begin{array}{c}
m \\l
\end{array}
\right ]_{q^2}\xi_2^{m-l}\La_s^{-m+l}\p_s^{m-l}\xi_1^l\La_s^{-l}\p_s^l=
$$
$$
=\sum_{m=0}^\infty\frac{(1-q^2)^m}{(q^2,q^2)_m}
\sum_{l=0}^m
\left[\begin{array}{c}
m \\l
\end{array}
\right ]_{q^2}\xi_2^{m-l}\xi_1^l\La_s^{-m}\p_s^m=
$$
$$
=\sum_{m=0}^\infty\frac{(1-q^2)^m}{(q^2,q^2)_m}
(\xi_1+\xi_2)^m\La_s^{-m}\p_z^m=
e_{q^2}((1-q^2)(\xi_1+\xi_2)\La_s^{-1}\p_z).
$$
\rule{5pt}{5pt}

\noindent{\bf Rule (O)}(Order). {\it If a function depends on several
variables then it is necessary to put them in order according to
subordination before one takes its restriction on the lattice
$\{q^{2n}\}$ so that the subordinate variable stands to the right.}
\begin{predl}\label{p4.3}
Operator $T_\xi$ for ${\xi=q^{2m}}$ can be represented by the form
\begin{equation}\label{4.5}
T_{q^{2m}}=\sum_{k=0}^\infty\frac{q^{2k(k+m)}}{(q^2,q^2)_k}s^{-k}
E_{q^2}(-q^{2(m+k+1)}s^{-1})\La_s^{-k}.
\end{equation}
\end{predl}
{\bf Proof}. It is easily to prove by the induction on $k$ that for an
arbitrary $k\ge0$
$$
\p_s^k\psi(s)=\frac1{(1-q^2)^ks^k}\sum_{l=0}^k(-1)^lq^{-l(2k-l-1)}
\left[\begin{array}{c}
k \\l
\end{array}
\right]_{q^2}\psi(q^{2l}s).
$$
It follows from (\ref{4.1}) and (\ref{4.2}) that
$$
T_\xi\psi(s)=\sum_{k=0}^\infty\xi^kq^{2k^2}s^{-k}
\sum_{l=0}^k(-1)^l\frac{q^{-l(2k-l-1)}}{(q^2,q^2)_l(q^2,q^2)_{k-l}}
\psi(q^{2l-2k}s)=
$$
$$
=\sum_{l=0}^\infty\frac{(-1)^lq^{l(l+1)}}{(q^2,q^2)_l}
\sum_{k=0}^\infty\frac{q^{2k(k+l)}}{(q^2,q^2)_k}\xi^{k+l}
s^{-k-l}\psi(q^{-2k}s)=
$$
$$
\sum_{k=0}^\infty\frac{q^{2k^2}}{(q^2,q^2)_k}\xi^k
\left(\sum_{l=0}^\infty(-1)^l\frac{q^{l(l+2k+1)}}{(q^2,q^2)_k}
\xi^lz^{-l}\right)s^{-k}\psi(q^{-2k}s).
$$
It is possible to change the order of summation because the inner series
are converged uniformly with respect to $k$ and $l$. Substituting
$\xi=q^{2m}$ we obtain
$$
T_{q^{2m}}\psi(s)=
\sum_{k=0}^\infty\frac{q^{2k(m+k)}s^{-k}}{(q^2,q^2)_l}\psi(q^{-2k}s)
\sum_{l=0}^\infty(-1)^l\frac{q^{l(l-1)}q^{2l(m+k+1)}s^{-l}}
{(q^2,q^2)_l}=
$$
\begin{equation}\label{4.6}
=\sum_{k=0}^\infty\frac{q^{2k(m+k)}}{(q^2,q^2)_l}s^{-k}
E_{q^2}(-q^{2(m+k+1)}s^{-1})\psi(q^{-2k}s).
\end{equation}
\rule{5pt}{5pt}

\begin{cor}\label{c4.1}
If $\psi(s)$ is limited by a constant $C$ i.e. for an arbitrary $n
~~~|\psi(q^{2n})|\le C$ then for any $n$ and $m$
$$
|T_{q^{2m}}\psi(q^{2n})|\le C.
$$
\end{cor}
{\bf Proof}. It is seen from (\ref{2.2}) that if $k>0$ then
$E_{q^2}(-q^{2k})>0$, and if $k\le0$ then $E_{q^2}(-q^{2k})=0$.
Substituting $s=q^{2n}$ in (\ref{4.5}) we obtain
$$
|T_{q^{2m}}\psi(q^{2n})|\le
C\sum_{k=0}^\infty\frac{q^{2k(k+m-n)}}{(q^2,q^2)_k}
E_{q^2}(q^{2(k+m+1-n)})=
$$
$$
=C\sum_{k=0}^\infty\frac{q^{2k(k+m-n)}}{(q^2,q^2)_k}
\sum_{l=0}^\infty\frac{(-1)^lq^{l(l+1)}q^{2l(k+m-n)}}{(q^2,q^2)_l}=
C\sum_{k=0}^\infty\frac{q^{2k(k+m-n)}}{(q^2,q^2)_k}
\sum_{l=0}^k(-1)^l\left[\begin{array}{c}
k \\l
\end{array}
\right ]_{q^2}q^{-l(2k-l-1)}.
$$
By induction on $k$
$$
\sum_{l=0}^k(-1)^l\left[\begin{array}{c}
k \\l
\end{array}
\right ]_{q^2}q^{-l(2k-l-1)}=\left\{\begin{array}{rcl}
1& {\rm for}&k=0\\ 0&{\rm for}&k\ge1.\\
\end{array}
\right.
$$
From here the statement of the Corollary follows.\rule{5pt}{5pt}

\begin{predl}\label{p4.4}
If $\psi(s)\in S^{q^2}$ then
\begin{equation}\label{4.7}
\int T_\xi\psi(s)d_{q^2}s=\int\psi(s)d_{q^2}s.
\end{equation}
\end{predl}
{\bf Proof.} It follows from (\ref{4.2}) that
$$
T_\xi\psi(s)=\psi(s)+\sum_{k=1}^\infty\frac{(1-q^2)^k}
{(q^2,q^2)_k}\xi^k\La_s^{-k}\p_s^k\psi(s).
$$
(\ref{4.7}) follows from Proposition \ref{p3.1}.
\rule{5pt}{5pt}

\begin{predl}\label{p4.5}
For $\psi(s)\in S^{q^2}$, for any $n\ge0,~~ m\ge0$ and for any $r,~~t$
$$
|s^n\p_s^mT_\xi\psi(s)|_{s=q^{2r},\xi=q^{2t}}|\le C_{n,m}.
$$
\end{predl}
{\bf Proof}. By induction on $k$
$$
\p_s^k(s\psi(s))=\frac{1-q^{2k}}{1-q^2}\p_s^{k-1}\psi(s)+
q^{2k}s\p_s^k\psi(s).
$$
We get from here and (\ref{4.2})
$$
sT_\xi\psi(s)=T_{q^{-2}\xi}(s\psi(s))-
q^{-2}\xi\La_s^{-1}T_{q^{-2}\xi}\psi(s).
$$
By induction on $n$
\begin{equation}\label{4.8}
s^nT_\xi\psi(s)=\sum_{l=0}^n(-1)^lC_n^lq^{-2l}
\xi^l\La_s^{-l}T_{q^{-2n}\xi}(\psi(s)s^{n-l}).
\end{equation}
It follows from (\ref{4.1}) that for any $m\ge0$
\begin{equation}\label{4.9}
\p_s^mT_\xi\psi(s)=T_\xi\p_s^m\psi(s).
\end{equation}
Proposition follows from (\ref{4.8}), (\ref{4.9})
and from Corollary \ref{c4.1}.\rule{5pt}{5pt}

The next theorem follows from Proposition \ref{p4.5}
\begin{theo}\label{t4.1}
The translation operator $T_\xi$ is the bounded operator in the space
$S^{q^2}$.
\end{theo}

Using (\ref{3.4}), (\ref{3.5}) we define the conjugate
operator $T_\xi^*$ in the space of the $q^2$-distributions $(S^{q^2})'$:
\begin{equation}\label{4.10}
T_\xi^*=e_{q^2}(-(1-q^2)q^2\xi\p_s), ~~~~~\xi s=q^2s\xi.
\end{equation}

If $\xi_1\xi_2=q^2\xi_2\xi_1$ then
$$
T_{\xi_2}^*T_{\xi_1}^*=T_{\xi_1+\xi_2}^*.
$$

\section{$q^2$-convolution}
\setcounter{equation}{0}

Let $g(s)$ be the $q^2$-distribution determined on the space $S^{q^2}$
of functions of one variable $s$, and $r(\xi)$ be the $q^2$-distribution
determined on the space $S^{q^2}$ of functions of one variable $\xi$,
moreover $s$ and $\xi$ are connected by relation (\ref{4.1}). We keep
the designation $S^{q^2}$ for the space of functions of two variables
$\psi(\xi,s)$. Then the functional
$$
r(\xi)\times g(s)
$$
is well-defined on this space and we call it the direct product of the
functionals $r(\xi)$ and $g(s)$
$$
<g(s),<r(\xi),\psi(\xi,s)>>=<r(\xi),<g(s),\psi(\xi,s)>>.
$$
In addition if the functionals $r$ and $g$ are regular then we must
succeed to Rule ({\bf O}).

\begin{defi}\label{d5.1}
We will call the functional
\begin{equation}\label{5.1}
<r*g,\psi>=<r(\xi)\times g(s),T_\xi\psi(s)>
\end{equation}
by the $q^2$-convolution of two $q^2$-distributions from $(S^{q^2})'$.
\end{defi}

\begin{defi}\label{d5.2}
$h(z)$ is the multiplicator in $S_{q^2}$ if for an arbitrary
function $\phi(z)\in S_{q^2} ~~~h(z)\phi(z)\in S_{q^2}$.
\end{defi}
\begin{predl}\label{p5.1}
The multiplication on an infinitely $q^2$-differentiable function $h(z)$
complying with inequality
\begin{equation}\label{5.2}
|\p_z^kh(z)|\le C_k(1+|z|^l)
\end{equation}
for some $l\ge0$ is a bounded operator in $S_{q^2}$.
\end{predl}
{\bf Proof}. By induction on $m$
$$
\p_z^m[h(z)\phi(z)]=\sum_{i=0}^m\left[
\begin{array}{c} m \\i
\end{array}
\right ]_{q^2}q^{-2i(m-i)}\p_z^{m-i}h(q^{2i}z)\p_z^i\phi(z).
$$
Then for any $n\ge0,~~~m\ge0$
$$
|z^n\p_z^m[h(z)\phi(z)]|=|z^n\sum_{i=0}^m\left[
\begin{array}{c} m \\i
\end{array}
\right ]_{q^2}q^{-2i(m-i)}\p_z^{m-i}h(q^{2i}z)\p_z^i\phi(z)|\le
$$
$$
\le\sum_{i=0}^m\left[
\begin{array}{c} m \\i
\end{array}
\right ]_{q^2}q^{-2i(m-i)}|\p_z^{m-i}h(q^{2i}z)||z^m\p_z^i\phi(z)|.
$$
It follows from (\ref{3.1}) and (\ref{5.2}) that
$$
|z^n\p_z^m[u(z)\phi(z)]|\le\sum_{i=0}^m\left[
\begin{array}{c} m \\i
\end{array}
\right ]_{q^2}q^{-2i(m-i)}C_{m-i}(1+|z|^l)|z^m\p_z^i\phi(z)|\le
$$
$$
\le\sum_{i=0}^m\left[
\begin{array}{c} m \\i
\end{array}
\right ]_{q^2}q^{-2i(m-i)}C_{m-i}C_{m,i}+
\sum_{i=0}^m\left[
\begin{array}{c} m \\i
\end{array}
\right ]_{q^2}q^{-2i(m-i)}C_{m-i}C_{m+l,i}=\tilde C_{m+l,n}.
$$
\rule{5pt}{5pt}

Thus the functions satisfying (\ref{5.2}), i.e. growing no
faster then some power of $|z|$ for $|z|\to\infty$, are the
multiplicator in $S_{q^2}$.

\begin{defi}\label{d5.3}
We will call functional $r\in (S^{q^2})'$ by $q^2$-convolutor in
$S^{q^2}$, if for an arbitrary function $\psi(s)\in S^{q^2}$ the
$q^2$-convolution
$$
(r*\psi)(s)=\int
d_{q^2}\xi\overline{r(\xi)}T_\xi\psi(s)
$$
exists and belongs to
$S^{q^2}$.
\end{defi}
\begin{predl}\label{p5.2}
If functional $h\in S_{q^2}$' is the multiplicator in $S_{q^2}$, then
its $q^2$-Fourier transform ${\cal F}_{q^2}'h=r\in(S^{q^2})'$ is the
$q^2$-convolutor in $S^{q^2}$.
\end{predl}
{\bf Proof.} Let functional $h$ correspond to the infinitely
$q^2$-differentiable function $h(z)$, and $zh(z)$ be an absolutely
$q^2$-integrable. Then functional ${\cal F}_{q^2}'h =r$ corresponds to the
function $r(s)$ and for any $k$ and $n\ge0$ a constant $C_n>0$ exists so
that $q^{2kn}|r(\pm q^{2k})|\le C_n$.  It follows from here and from
Theorem \ref{t4.1} that the series
$$
(1-q^2)\sum_{k=-\infty}^\infty q^{2k}[\overline{r(q^{2k})}
T_{q^{2k}}\psi(s)+\overline{r(-q^{2k})}T_{-q^{2k}}\psi(s)]
$$
converges uniformly with respect to $s$ together with $q^2$-derivatives
with respect to $s$ and we can $q^2$-differentiate it term-by-term.

It is easily to show as in the proof of Proposition \ref{p4.5} that for
any $n\ge0$ and $m\ge0$
\begin{equation}\label{5.3}
s^n\p_s^mT_{\pm q^{2k}}\psi(s)=\sum_{l=0}^n(-1)^l\left[
\begin{array}{c} n \\l
\end{array}
\right ]_{q^2}q^{2l(k-m)}\La_s^{-l}T_{\pm q^{2(k-m)}}(s^{n-l}\psi(s)).
\end{equation}
Hence, for any $n\ge0$ and $m\ge0$ and for an arbitrary
$\psi(s)\in S^{q^2}$
$$
|s^n\p_s^m(r*\psi)(s)|\le
(1-q^2)\sum_{k=-\infty}^\infty q^{2k}[|\overline{r(q^{2k})}
s^n\p_s^mT_{q^{2k}}\psi(s)|+|\overline{r(-q^{2k})}
s^n\p_s^mT_{-q^{2k}}\psi(s)|].
$$
Proposition follows from Corollary \ref{c4.1} and (\ref{5.3}).

Let now functional $h$ correspond to infinitely $q$-differentiable
function $h(z)$ growing for $|z|\to\infty$ no faster then $|z|^p, ~~p>0$
is integer . That is $h(z)$ has form
$$
h(z)=\sum_{k=0}^pz^kh_k(z),
$$
where $h_k(z)$ are such infinitely $q^2$-differentiable functions that
$zh_k(z)$ are absolutely $q^2$-integrable. It follows from (\ref{3.12})
that
$$
r(s)=\sum_{k=0}^p(-i\p_s)^kr_k(s),~~~r_k(s)={\cal F}_{q^2}'h_k(z).
$$
On the other hand using the formula of $q^2$-integration by parts
(\ref{3.2}) it is easily to obtain
\begin{equation}\label{5.4}
\int d_{q^2}\xi\overline{\p_\xi^kr_k(\xi)}T_\xi\psi(s)=(-1)^kq^{k(k-1)}
\int d_{q^2}\overline{r_k(s)}T_\xi(\p_s^k\psi(q^{-2k}s)).
\end{equation}
So the statement of the Proposition is truly in this case also.
\rule{5pt}{5pt}

The formula
\begin{equation}\label{5.5}
(\p_xi^kr*\psi)(s)=(-1)^{-k(k+1)}\La_s^{-k}\p_s^k(r*\psi)(s).
\end{equation}
follows from (\ref{5.4})

\begin{predl}\label{p5.3}
If $r$ is the $q^2$-convolutor in $S^{q^2}$, then the $q^2$-convolution
$(r*\psi)(s)$ is the $q^2$-Fourier transform of the product
$\overline{h(z)}\phi(z)$, where $\phi(z)={\cal F}_{q^2}^{-1}\psi(s)$ and
$h(z)=({\cal F}_{q^2}')^{-1}r(s)$.
\end{predl}
{\bf Proof.} It was proved in \cite{OR} (Lemma (5.2)) that
\begin{equation}\label{5.6}
\int E_{q^2}(-i(1-q^2)zs)E_{q^2}(i(1-q^2)q^2s)d_{q^2}s=\left\{
\begin{array}{rcl}
\frac2{1-q^2}\Theta_0 &{\rm for}&z=1\\0&{\rm for}&z\ne1.\\
\end{array}\right.
\end{equation}

It is easily to convince that
\begin{equation}\label{5.7}
e_{q^2}((1-q^2)\xi\La_s^{-1}\p_s)E_{q^2}(i(1-q^2)q^{2n+2}s)=
E_{q^2}(i(1-q^2)q^{2n+2}\xi)E_{q^2}(i(1-q^2)q^{2n+2}s).
\end{equation}

At first let $r$ be a regular functional corresponding to function $r(s)$.
It follows from (\ref{3.7}) and (\ref{3.10})
$$
\int d_{q^2}\xi\overline{r(\xi)}e_{q^2}((1-q^2)\xi\La_s^{-1}\p_s)\psi(s)=
$$
$$
=\frac{(1-q^2)^2}{2\Theta_0}\int d_{q^2}\xi\sum_{m=-\infty}^\infty q^{2m}
[\overline{h(q^{2m})}e_{q^2}(-i(1-q^2)q^{2m}\xi)+
\overline{h(-q^{2m})}e_{q^2}(i(1-q^2)q^{2m}\xi)]\times
$$
$$
\times e_{q^2}((1-q^2)\xi\La_s^{-1}\p_s)\sum_{n=-\infty}^\infty q^{2n}
[\phi(q^{2n})E_{q^2}(i(1-q^2)q^{2n+2}s)+
\phi(-q^{2n})E_{q^2}(-i(1-q^2)q^{2n+2}s)].
$$
These series converge uniformly with respect to $\xi$ and so we can
$q^2$-integrate them term-by-term. Then using (\ref{5.6}) and (\ref{5.7}),
we obtain
$$
\int d_{q^2}\xi\overline{r(\xi)}e_{q^2}((1-q^2)\xi\La_s^{-1}\p_s)\psi(s)=
\frac{(1-q^2)^2}{2\Theta_0}\sum_{m=-\infty}^\infty q^{2m}
\sum_{n=-\infty}^\infty q^{2n}\times
$$
$$
\times[\overline{h(q^{2m})}\phi(q^{2n})\int d_{q^2}\xi
e_{q^2}(-i(1-q^2)q^{2m}\xi)
E_{q^2}(i(1-q^2)q^{2n+2}\xi)E_{q^2}(i(1-q^2)q^{2n+2}s)+
$$
$$
+\overline{h(q^{2m})}\phi(-q^{2n})\int d_{q^2}\xi
e_{q^2}(-i(1-q^2)q^{2m}\xi)
E_{q^2}(-i(1-q^2)q^{2n+2}\xi)E_{q^2}(-i(1-q^2)q^{2n+2}s)+
$$
$$
+\overline{h(-q^{2m})}\phi(q^{2n})\int d_{q^2}\xi
e_{q^2}(i(1-q^2)q^{2m}\xi)
E_{q^2}(i(1-q^2)q^{2n+2}\xi)E_{q^2}(i(1-q^2)q^{2n+2}s)+
$$
$$
+\overline{h(-q^{2m})}\phi(-q^{2n})\int d_{q^2}\xi
e_{q^2}(i(1-q^2)q^{2m}\xi)
E_{q^2}(-i(1-q^2)q^{2n+2}\xi)E_{q^2}(-i(1-q^2)q^{2n+2}s)]=
$$
$$
=(1-q^2)\sum_{n=-\infty}^\infty q^{2n}\left[
\overline{h(q^{2n})}\phi(q^{2n})E_{q^2}(i(1-q^2)q^{2n+2}s)+
\overline{h(-q^{2n})}\phi(-q^{2n})E_{q^2}(-i(1-q^2)q^{2n+2}s)\right]=
$$
$$
=\int d_{q^2}z\overline{h(z)}\phi(z)
\phantom._0\Phi_1(-;0;q^2(1-q^2)q^2zs).
$$

Now let functional $r\in (S^{q^2})'$ be  a $q^2$-singular one with
the multiplicity $p$ (see Definition \ref{d3.2}), i.e.
$$
r(s)=\sum_{k=0}^p\p_s^kr_k(s),
$$
where $r_k(s)$ are regular functionals. It follows from (\ref{3.12}) that
\begin{equation}\label{5.8}
({\cal F}_{q^2}')^{-1}r=h(z)=\sum_{k=0}^p(iz)^kh_k(z).
\end{equation}
We have from (\ref{3.7}) and (\ref{5.5})
$$
{\cal F}_{q^2}^{-1}(r*\psi)(s)=
{\cal F}_{q^2}^{-1}\sum_{k=0}^p(-1)^kq^{-k(k+1)}
\La_s^{-k}\p_s^k(r_k*\psi)(s)=
$$
$$
={\cal F}_{q^2}^{-1}\sum_{k=0}^p(-1)^k(q^{-2}\La_s^{-1}\p_s)^k(r_k*\psi)(s)=
=\sum_{k=0}^p(iz)^k\overline{h_k(z)}\phi(z).
$$
Then we obtain from (\ref{5.8})
$$
{\cal F}_{q^2}^{-1}(r*\psi)(s)=\overline{h(z)}\phi(z).
$$
\rule{5pt}{5pt}

\begin{theo}\label{t5.1}
Let $g\in(S^{q^2})'$ and $r(s)$ be a $q^2$-convolutor in $S^{q^2}$.
Then the $q^2$-convolution  $r*g$ (\ref{5.1}) is the $q^2$-Fourier
transform of the product of the $q^2$-distributions $hf$, where
$h=({\cal F}_{q^2}')^{-1}r$, and $f=({\cal F}_{q^2}')^{-1}g$.
\end{theo}
{\bf Proof.} For an arbitrary $\psi(s)$ and $\phi(z)=
{\cal F}_{q^2}^{-1}\phi(z)$ we have from the definition of the
$q^2$-distributions (\ref{5.1} and from the definition of the $q^2$-Fourier
transform of the $q^2$-distributions \ref{3.8}
\begin{equation}\label{5.9}
<r*g,\psi>=<g,r*\psi>=<f,\overline{h}\phi>=<fh,\phi>.
\end{equation}
\rule{5pt}{5pt}
\begin{cor}\label{c5.1}
If $g(s)$ and $r(s)$ are $q^2$-convolutors in $S^{q^2}$, then them
$q^2$-convolution is commutative one i.e. for an arbitrary $\psi(s)\in
S^{q^2}$
\begin{equation}\label{5.10}
<r*g,\psi>=<g*r,\psi>.
\end{equation}
\end{cor}
{\bf The proof} follows from commutativity of the product
$f(z)h(z)$ in right side of (\ref{5.9}). \rule{5pt}{5pt}

\begin{cor}\label{c5.2}
If $g(s)$ is $q^2$-convolutor in $S^{q^2}$, then
\begin{equation}\label{5.11}
g*\delta_{q^2}=\delta_{q^2}*g=g.
\end{equation}
\end{cor}
{\bf Proof}. $q^2-\delta$-function is the $q^2$-convolutor
in $S^{q^2}$ because it is the $q^2$-Fourier transform of
$\frac12$, which is the multiplicator in $S_{q^2}$ (see \cite{OR}). It
is follows from (\ref{4.10})
\begin{equation}\label{5.12}
g*\delta_{q^2}=\int d_{q^2}\xi\delta_{q^2}(\xi)T_\xi^*g(s)=T_0^*g(s)=g(s).
\end{equation}
\rule{5pt}{5pt}

\begin{predl}\label{p5.4}
\begin{equation}\label{5.13}
\p_sg*r=g*\p_sr.
\end{equation}
\end{predl}
{\bf Proof}. It is follows from (\ref{4.10})
$$
\p_\xi e_{q^2}(-(1-q^2)q^2\xi\p_s)=-q^2e_{q^2}(-(1-q^2)q^4\xi\p_s)\p_s.
$$
Now (\ref{5.13}) follows from (\ref{3.2}) and (\ref{3.3}).\rule{5pt}{5pt}

\section{The $q^2$-pseudo differential operators}
\setcounter{equation}{0}

Consider the $q^2$-distribution $s_+^{\nu-1}$ \cite{OR}. For an arbitrary
fixed $n\ge0$
$$
\int_0^\infty s^{\nu-1}\psi(s)d_{q^2}s=\int_0^1s^{\nu-1}\Bigl[\psi(s)-
\sum_{k=0}^n\frac{s^k}{k!}\psi^{(k)}(0)\Bigr]d_{q^2}s+
$$
$$
+\int_1^\infty s^{\nu-1}\psi(s)d_{q^2}s+(1-q^2)\sum_{k=0}^n
\frac1{k!(1-q^{2(\nu+k)})}\psi^{(k)}(0).
$$
It is seen from (\ref{3.2}), that the last sum can be represented by the
form
$$
\sum_{k=0}^n\frac{(1-q^2)^{k+1}}
{(q^2,q^2)_k(1-q^{2(\nu+k)})}\p_s^k\psi(0).
$$
On the other hand
$$
<\delta_{q^2},\p_s^k\psi>=(-1)^kq^{k(k+1)}<\p_s^k\delta_{q^2},\psi>.
$$
Hence $s_+^{\nu-1}$ is the meromorphic function of $\nu$ with the ordinary
poles  $\nu=-k, ~~k=0, 1,\ldots$ with the residues
\begin{equation}\label{6.1}
\res_{\nu=-k}s_+^{\nu-1}=(-1)^kq^{k(k+1)}\frac{(1-q^2)^{k+1}}
{(q^2,q^2)_k}\p_s^k\delta_{q^2}(s).
\end{equation}

Consider now the $q^2$-$\G$-function
\begin{equation}\label{6.2}
\G_{q^2}(\nu)=\frac{(q^2,q^2)_\infty}{(q^{2\nu},q^2)_\infty}(1-q^2)^{1-\nu}.
\end{equation}
Obviously it is the meromorphic function with the ordinary poles $\nu=-k, ~~
k=0, 1,\ldots$ with the residues
$$
\res_{\nu=-k}\G_{q^2}(\nu)=
\lim_{\nu\to-k}\frac{(q^2,q^2)_\infty(1-q^2)^{1-\nu}}
{(1-q^{2\nu})\ldots(1-q^{2\nu+2k-2})(q^{2\nu+2k+2},q^2)_\infty}=
$$
\begin{equation}\label{6.3}
=\frac{(1-q^2)^{k+1}}{(1-q^{-2\nu})\ldots(1-q^{-2})}=
(-1)^kq^{k(k+1)}\frac{(1-q^2)^{k+1}}{(q^2,q^2)_k}.
\end{equation}
The next Proposition follows from (\ref{6.1}) and (\ref{6.3}).
\begin{predl}\label{p6.1}
The $q^2$-distribution $\frac{s_+^{\nu-1}}{\G_{q^2}(\nu)}$ is the entire
function of $\nu$ and
\begin{equation}\label{6.4}
\frac{s_+^{\nu-1}}{\G_{q^2}(\nu)}\Bigr|_{\nu=-k}=
\p_s^k\delta_{q^2}(s), ~~~k=0,\ldots.
\end{equation}
\end{predl}
\begin{defi}\label{d6.1}
For an arbitrary $q^2$-distribution $g\in (S^{q^2})'$ concentrated on the
lattice $\{q^{2n}\}$ and absolutely $q^2$-integrable on any segment we
will call the $q^2$-convolution
\begin{equation}\label{6.5}
g*\frac{s_+^{\nu-1}}{\G_{q^2}(\nu)}
\end{equation}
by the $q^2$-derivative of $g$ $-\nu$ order if $\nu<0$ and by the
$q^2$-primitive of $g$ $\nu$ order if $\nu>0$.
\end{defi}
In accordance with Definition \ref{d6.1} we introduce the
designation
$$
\p_s^{-\nu}g=g*\frac{s_+^{\nu-1}}{\G_{q^2}(\nu)}.
$$
Moreover the $q^2$-derivative $-\nu$ order is the $q^2$-primitive $\nu$
order if $\nu>0$.

For the proof of correctness of Definition \ref{d6.1} it is sufficient
to check its correctness for $\nu=0, \pm1$ and to prove the following
Proposition
\begin{predl}\label{p6.2}
For arbitrary $\nu$ and $\mu$
\begin{equation}\label{6.6}
\frac{s_+^{\nu-1}}{\G_{q^2}(\nu)}*\frac{s_+^{\mu-1}}{\G_{q^2}(\mu)}=
\frac{s_+^{\nu+\mu-1}}{\G_{q^2}(\nu+\mu)}.
\end{equation}
\end{predl}

Let $\nu=0$. It follows from (\ref{6.4}) and (\ref{5.11})
$$
g(s)*\frac{s_+^{\nu-1}}{\G_{q^2}(\nu)}\Bigr|_{\nu=0}=(g*\delta_{q^2})(s)=
g(s).
$$
Let $\nu=-1$. It follows from (\ref{6.4}), (\ref{5.11}) and (\ref{5.13})
$$
g(s)*\frac{s_+^{\nu-1}}{\G_{q^2}(\nu)}\Bigr|_{\nu=-1}=
(g*\p_s\delta_{q^2})(s)=(\delta_{q^2}*\p_sg)(s)=\p_sg(s).
$$
Finally let $\nu=1$. As $g$ is concentrated on the lattice $\{q^{2n}\}$
and
$$
\Bigl(T_\xi^*s^0\Bigr)_+=\left\{
\begin{array}{rcl}
1 & {\rm }& \xi\le s \\
0 & {\rm }& \xi>s,\\
\end{array}
\right.
$$
then
$$
g(s)*\frac{s_+^0}{\G_{q^2}(1)}=\int_0^sd_{q^2}\xi g(\xi)=
(1-q^2)s\sum_{m=0}^\infty q^{2m}g(q^{2m}s).
$$
So the last function is the $q^2$-primitive of $g(s)$ because its
$q^2$-derivative is $g(s)$.

{\bf Proof} of Proposition \ref{p6.2}. The statement of Proposition is
trivial if $\nu$ and $\mu$ are integer. Let $\nu$ is not integer.
Obviously
$$
\p_s^ks^{\nu-1}=(-1)^kq^{k(2\nu-k-1)}\frac{(q^{-2\nu+2},q^2)_k}
{(1-q^2)^k}s^{\nu-k-1}.
$$
Using the properties of the $q^2$-binomial formula \cite{GR}
$$
\sum_{k=0}^\infty\frac{(a,q^2)_k}{(q^2,q^2)_k}x^k=
\frac{(ax,q^2)_\infty}{(x,q^2)_\infty},
$$
and (\ref{4.10}) for function $g(s)$ concentrated on the set $s\ge0$
we obtain
$$
\int_0^sd_{q^2}\xi g(\xi)T_\xi^*\frac{s^{\nu-1}}{\G_{q^2}(\nu)}=
\int_0^sd_{q^2}\xi g(\xi)\sum_{k=0}^\infty q^{2\nu k}
\frac{(q^{-2\nu+2},q^2)_k}{(q^2,q^2)_k}\frac{s^{\nu-k-1}}{\G_{q^2}(\nu)}=
$$
$$
=(1-q^2)\frac{s^\nu}{\G_{q^2}(\nu)}\sum_{m=0}^\infty q^{2m}g(q^{2m}s)
\sum_{k=0}^\infty\frac{(q^{-2\nu+2},q^2)_k}{(q^2,q^2)_k}q^{2k(\nu+m)}=
$$
$$
=(1-q^2)\frac{s^\nu}{\G_{q^2}(\nu)}\sum_{m=0}^\infty
\frac{q^{2m}(q^{2m+2},q^2)_\infty}{(q^{2(\nu+m},q^2)_\infty}g(q^{2m}s)=
$$
$$
=(1-q^2)\frac{s^\nu}{\G_{q^2}(\nu)}\frac{(q^2,q^2)_\infty}
{(q^{2\nu},q^2)_\infty}\sum_{m=0}^\infty
\frac{q^{2m}(q^{2\nu},q^2)_m}{(q^2,q^2)_m}g(q^{2m}s)=
(1-q^2)^\nu s^\nu\sum_{m=0}^\infty
\frac{q^{2m}(q^{2\nu},q^2)_m}{(q^2,q^2)_m}g(q^{2m}s).
$$
So
\begin{equation}\label{6.7}
\p_s^{-\nu}g(s)=(1-q^2)^\nu s^\nu\sum_{m=0}^\infty
\frac{q^{2m}(q^{2\nu},q^2)_m}{(q^2,q^2)_m}g(q^{2m}s).
\end{equation}

It is easily to prove by the induction on $k$ that for any non integer
$\nu$ and $\mu$ and for any integer $k\ge0$
\begin{equation}\label{6.8}
\sum_{l=0}^kq^{2\nu l} \left[\begin{array}{c} k \\l
\end{array}
\right]_{q^2}
(q^{2\mu},q^2)_l(q^{2\nu},q^2)_{k-l}=(q^{2\mu+2\nu},q^2)_k.
\end{equation}

Consider now the composition of the operators $\p_s^{-\nu}$ and
$\p_s^{-\mu}$ for non integer $\nu$ and $\mu$. It follows from (\ref{6.7})
$$
\p_s^{-\mu}\p_s^{-\nu}g(s)=(1-q^2)^\mu s^\mu\sum_{l=0}^\infty
\frac{q^{2l}(q^{2\mu},q^2)_l}{(q^2,q^2)_l}
(1-q^2)^\nu(q^{2l}s)^\nu\sum_{m=0}^\infty
\frac{q^{2m}(q^{2\nu},q^2)_m}{(q^2,q^2)_m}g(q^{2m+2l}s)=
$$
$$
=(1-q^2)^{\nu+\mu}s^{\nu+\mu}\sum_{k=0}^\infty
\frac{q^{2n}(q^{2\nu+2\mu},q^2)_k}{(q^2,q^2)_k}g(q^{2k}s)
\sum_{l=0}^kq^{2\nu l}
\left[\begin{array}{c}
k \\l
\end{array}
\right ]_{q^2}
\frac{(q^{2\mu},q^2)_l(q^{2\nu},q^2)_{k-l}}{(q^{2\mu+2\nu},q^2)_k}=
$$
$$
=(1-q^2)^{\nu+\mu}s^{\nu+\mu}\sum_{k=0}^\infty
\frac{q^{2n}(q^{2\nu+2\mu},q^2)_k}{(q^2,q^2)_k}g(q^{2k}s)=
\p_s^{-(\nu+\mu)}g(s),
$$
The last equality follows from (\ref{6.8}). Proposition is proved.
\rule{5pt}{5pt}

In \cite{OR} the $q^2$-Fourier transforms of the $q^2$-distributions from
the space $S_{q^2}'~~~z_+^\nu$ and $z_-^\nu$ are calculated. It is
follows from these formulas
\begin{equation}\label{6.9}
\frac{s_+^{\nu-1}}{\G_{q^2}(\nu)}
=2\Theta_0{\cal F}_{q^2}'[A_\nu z_+^{-\nu}+\bar A_\nu z_-^{-\nu}],
\end{equation}
where $\Theta_0$ is determined by (\ref{2.4}),
\begin{equation}\label{6.10}
A_\nu=\frac{c_\nu(1-q^2)^{\nu-1}}{c_\nu^2-\bar c_\nu^2},
\end{equation}
\begin{equation}\label{6.11}
c_\nu=\sum_{m=-\infty}^\infty\frac{q^{-2\nu m}+i(1-q^2)q^{2(1-\nu)m}}
{(1-q^2)^{-1}q^{-2m}+(1-q^2)q^{2m}}.
\end{equation}

The next theorem is follows from Theorem \ref{t5.1} and (\ref{6.9})
\begin{theo}\label{t6.1}(Addition theorem)
$$
A_\nu A_\mu=\frac1{2\Theta_0}A_{\nu+\mu}.
$$
\end{theo}

\bigskip
\small{

}

\end{document}